\documentclass[10pt]{article}
\usepackage{amssymb,amsfonts, amsmath}

\newtheorem{theorem}{Theorem}
\newtheorem{lemma}[theorem]{Lemma}

\newtheorem{corollary}[theorem]{Corollary}

\newcommand{\newsection}[1] {\section{#1}\setcounter{theorem}{0} 
 \setcounter{equation}{0}\par\noindent}

\newcommand{\R}{{\mathbb R}}
\newcommand{\I}{{\rm I}}
\newcommand{\grad}{\nabla}

\renewcommand{\l}{\lambda}
\renewcommand{\div}{\text{div}}
\newcommand{\bc}{{\mathbf c}}
\newcommand{\bs}{{\mathbf s}}
\newcommand{\br}{{\mathbf r}}
\newcommand{\rc}{{\mathrm c}}
\newcommand{\rs}{{\mathrm s}}

\title
{Decoupling of modes for the elastic wave equation\\
in media of limited smoothness}

\author{Valeriy Brytik$^1$ \and Maarten V. de Hoop$^1$ \and Hart F. Smith$^2$ \and Gunther Uhlmann$^2$}

\begin{document}

\maketitle

\noindent
{${}^1$Center for Computational and Applied Mathemematics, and Geo-Mathematical Imaging Group, Purdue University, West Lafayette, Indiana, USA.}

\medskip

\noindent
{${}^2$Department of Mathematics, University of Washington, Seattle, Washington, USA.}

\bigskip

\begin{abstract}
\noindent
We establish a decoupling result for the $P$ and $S$ waves of linear, isotropic
elasticity, in the setting of twice-differentiable Lam\'e parameters. Precisely,
we show that the $P\leftrightarrow S$ components of the wave propagation operator
are regularizing of order one on $L^2$ data, by establishing the
diagonalization of the elastic system modulo a $L^2$-bounded operator.
Effecting the diagonalization in the setting of twice-differentiable coefficients
depends upon the symbol of the conjugation operator having a particular structure.
\end{abstract}

\newsection{Introduction}
We consider the linear, isotropic elastic wave equation on $\R^n$, and Cauchy
initial value problem on a time interval $[-t_0,t_0]$,
\begin{equation} \label{eq:CIVP}
\left\{\begin{array}{rcll}
   \partial_t^2w_i-(A(x,D) w)_i &=& 0\,,
\\[0.4cm]
   w_i |_{t=0} &=& f_i\,,
\\[0.4cm]
   \partial_t w_i |_{t=0} &=& g_i\,,
\end{array}\right.
\end{equation}
where $w=(w_1,\ldots,w_n)$ denotes the displacement vector,
and
\begin{equation}\label{Aform}
   (A(x,D)w)_i=\sum_{k=1}^n
   \partial_k\bigl(\,\mu(x)(\partial_kw_i+\partial_iw_k)\,\bigr)
  +\partial_i\bigl(\,\l(x)\partial_kw_k\bigr)\,.
\end{equation}
The elastic wave equation arises in many physical situations,
including the propagation of seismic waves through the Earth. The
identification and separation of $P$ and $S$ polarized constituents
(so-called phases) underlies many applications of seismic data
processing and imaging, which traditionally has been based on the
assumption of smoothness of the medium and asymptotic
considerations. As an example, we mention the notion of receiver
functions, which were introduced and developed by Vinnik (1977)
and Langston (1979).

Throughout this paper, we assume that the Lam\'e parameter functions 
$\mu$ and $\l$ belong to $C^{1,1}(\R^n)$, 
in that the first derivatives $\nabla\mu$ and $\nabla\l$ are globally bounded, Lipschitz functions
on $\R^n$.

The principal symbol $\sigma_2(A)(x,\xi)$ of $A(x,D)$ takes the form
\begin{equation}\label{princsymbA}
\begin{split}
-\sigma_2(A)&=\mu(x)\,|\xi|^2\,\I+\bigl(\mu(x)+\l(x)\bigr)\,\xi\otimes\xi
\\
\rule{0pt}{14pt}
&=\,\bigr(2\mu(x)+\l(x)\bigr)\,|\xi|^2\,\sigma_0(\Pi_P)
+\mu(x)\,|\xi|^2\,\sigma_0(\Pi_S)
\end{split}
\end{equation}
where $\Pi_P$ and $\Pi_S$ are the projection operators
$$
\Pi_P w = \Delta^{-1}\grad\bigl(\div\, w)\,,\qquad
\Pi_S w = (\I-\Pi_P)w
$$
corresponding to the Fourier multipliers
$$
\sigma_0(\Pi_P)_{il}=\frac{\xi_i\xi_l}{|\xi|^2}\,,\qquad
\sigma_0(\Pi_S)_{il}=\frac{|\xi|^2\delta_{il}-\xi_i\xi_l}{|\xi|^2}\,.
$$
The symbols of $\Pi_P$ and $\Pi_S$ are thus, for each $x,\xi$, projections
onto the eigenspaces of $A(x,\xi)$, respectively of dimension 1 and $n-1$.
This yields that
\begin{equation}\label{eq:princdiag}
\begin{split}
A(x,D)\,\Pi_Pw&=\bigl(2\mu(x)+\l(x)\bigr)\Delta \, \Pi_Pw+\text{l.o.t.}\,,\\
\rule{0pt}{16pt}
A(x,D)\,\Pi_Sw&=\mu(x)\Delta \, \Pi_Sw+\text{l.o.t.}\,,
\end{split}
\end{equation}
where $\text{l.o.t.}$ denotes an operator of first order in $w$ (on some finite
range of Sobolev spaces, as discussed in the next section.)

We call $\Pi_Pw$ and $\Pi_Sw$, respectively, the $P$ and $S$ modes of the elastic wave
$w$. If the Lam\'e parameters are constant, then 
the wave propagation operator commutes with $\Pi_P$ and $\Pi_S$,
so that wave evolution preserves modes which are of pure $P$ or pure
$S$ form.

In the case of nonconstant $\mu$ and $\l$, there will in general be
a nonvanishing coupling of the $P$ and $S$ modes in the evolution. 
The main result of this paper is that the $P\leftrightarrow S$
terms in the wave propagator are
one order smoother than the $P\rightarrow P$ and $S\rightarrow S$ terms.
Precisely, let the solution to \eqref{eq:CIVP}
be written in the form
$$
w(t,\cdot\,)=\bc(t)f+\bs(t)g\,,
$$
and decompose $\bc(t)=\bc_{PP}(t)+\bc_{PS}(t)+\bc_{SP}(t)+\bc_{SS}(t)$, 
where $\bc_{PP}(t)=\Pi_P\bc(t)\,\Pi_P$, etc. We similarly decompose $\bs(t)$.
The operators $\bc_{PP}(t)$ and $\bc_{SS}(t)$ are bounded on a range of Sobolev
spaces $H^s(\R^n)$, in particular for $0\le s\le 1$, and $\bs_{PP}(t)$, $\bs_{SS}(t)$
map $H^{s-1}(\R^n)$ to $H^s(\R^n)$ for the same range. (All maps are strongly
continuous in $t$.) We prove that
$$
\bc_{PS}(t)\,:\,L^2(\R^n)\rightarrow H^1(\R^n)\,, \qquad 
\bs_{PS}(t)\,:\,H^{-1}(\R^n)\rightarrow H^1(\R^n)\,.
$$
Additionally, we show that $\bc_{PP}(t)$ agrees with the propagator
for $A_{PP} = \Pi_P A \Pi_P$ up to a similarly regularizing term; see
Corollary \ref{cor:decoupled} for a precise statement.

Establishing regularization for the $P\leftrightarrow S$ components of
the propagator requires diagonalizing $A(x,D)$ to second order; that
is, one order beyond the principal symbol consideration of
\eqref{eq:princdiag}, since treating a first order derivative of $w$
as a driving force term will produce terms of the same order as $w$ in
the solution. This diagonalization is carried out through conjugation
by $I+K$, where $K$ is a compound pseudodifferential operator with
finite regularity symbol of order $-1$.

In the case of smooth $\mu$ and $\l$, a similar diagonalization was
carried out, in a more general setting, to all orders in Taylor (1975)
using the pseudodifferential calculus.  He used this
diagonalization to study the reflection of singularities for the
system of isotropic elasticity, as well as for other systems with
constant multiplicities. Stolk and de Hoop (2002)
have used Taylor's diagonalization in the study of the linearized seismic
inverse scattering problem. Dencker (1982)
analyzed the propagation
of polarization in Euclidean space for linear systems of principal
type with smooth coefficients, which includes isotropic
elasticity. Wang (1998)
applied Dencker's propagation of
polarization result to the linearized seismic inverse problem; in
particular he discusses the issue of polarization in reflections
caused by discontinuities of a perturbation of the medium.

In the case of $C^{1,1}$ coefficients, the construction of $K$ is more
delicate, and carrying through the composition calculus to first order
depends on the symbol of $K$ having a special form.  We believe that
this is the first article studying the interaction of $P$ and $S$
waves for the elastic system in non-smooth media, and it remains open
under what conditions this can be carried out in the $C^{1,1}$ setting
for more general diagonalizable systems with constant multiplicity
eigenvalues. We anticipate that our result will have applications in
the study of the linearized seismic inverse problem in media of
limited smoothness. There has been considerable activity in studying
unique continuation in the stationary case for non-smooth Lam\'e
parameters, see Lin et.~al.~(2010) 
and the references therein for the most
recent results.

\bigskip
\bigskip

\noindent{\bf \Large Notation.}

\medskip

All pseudodifferential operators that arise in this paper will have
symbols that are finite sums $p(x,\xi)=\sum_{j=1}^N p_{m-j}(x,\xi)$,
where $p_{m-j}$ is homogeneous in $\xi$ of degree $m-j$.  We use
$\sigma_m(p)$ to denote the principal symbol $p_m(x,\xi)$ of $p$,
where $m$ is the order of $p$.

We will use $L^2(\R^n)$ and $H^s(\R^n)$ to denote the $L^2$ and
$L^2$-based Sobolev spaces, in the case of complex valued functions as
well as $n$-dimensional vector valued functions, such as the solution
$w$.

Throughout the paper, $R$ will denote an operator that is bounded on
$L^2(\R^n)$; the form of $R$ may change in each occurrence.

\newsection{The Conjugation Operator}
In this section we construct the conjugation operator $K$, which diagonalizes
$A(x,D)$ up to a $L^2$-bounded operator.
We decompose $A=A_{PP}+A_{PS}+A_{SP}+A_{SS}$, where
\begin{align*}
A_{PP}=\Pi_P\,A(x,D)\,\Pi_P\,,\quad&A_{PS}=\Pi_P\,A(x,D)\,\Pi_S\,,\\
\rule{0pt}{16pt}
A_{SP}=\Pi_S\,A(x,D)\,\Pi_P\,,\quad&A_{SS}=\Pi_S\,A(x,D)\,\Pi_S\,.
\end{align*}
Note that $A_{SP}=A_{PS}^*$, while both $A_{SS}$ and $A_{PP}$ are self-adjoint.
Furthermore $A_{SP}$, and hence $A_{PS}$, is of order 1, in that
\begin{equation}\label{APSmap}
A_{PS}\,:\,H^{s+1}(\R^n)\rightarrow H^{s}(\R^n)\quad\text{for}\quad -2\le s\le 1\,.
\end{equation}
To see this, we write
\begin{equation}\label{APSform}
A_{PS}w=\Delta^{-1}\grad(\div \, A \,\Pi_Sw)=\Delta^{-1}\grad(\,[\div,A]\,\Pi_Sw)\,,
\end{equation}
where we use $\div \, \Pi_S\equiv 0$. Each component of the composition
then takes the form $\Delta^{-1}\grad(\partial(b\,\partial\,\Pi_S w))$,
where $b$ denotes a first-order derivative of $\mu$ or $\l$, which being
Lipschitz acts as a multiplier on $H^s$ for $-1\le s\le 1$.
This establishes \eqref{APSmap} for such $s$. If instead one writes
$$
A_{PS}w=-\Pi_P A\,(\Delta-\nabla\,\div)\,\Delta^{-1}w
$$
and commutes one derivative from the right past $A$,
one obtains \eqref{APSmap} for the range $-2\le s\le 0$.

Our main conjugation result is the following.
\begin{lemma}
There exists an operator $K$ of order $-1$, in the sense that
\begin{equation}\label{Kcond2}
K\,:\,H^{s}(\R^n)\rightarrow H^{s+1}(\R^n)\quad\text{for}\quad -1\le s\le 1\,,
\end{equation}
such that
\begin{align}
\Pi_P(I+K)A(x,D)w&=A_{PP}\,(I+K)w+R_Pw\,,\label{conj1}\\
\rule{0pt}{16pt}
\Pi_S(I+K)A(x,D)w&=A_{SS}\,(I+K)w+R_Sw\,,\label{conj1'}
\end{align}
where $R_P$ and $R_S$ are bounded maps on
$L^2(\R^n)$.

Additionally,
\begin{equation}\label{Kcond1}
K_{PP}=K_{SS}=0\,,\qquad K_{SP}^*=-K_{PS}\,,
\end{equation}
where $K_{PP}=\Pi_P\,K\,\Pi_P$, $\;K_{SP}=\Pi_S\,K\,\Pi_P$, etc.
\end{lemma}

To establish \eqref{conj1}-\eqref{conj1'}, we expand $A$ and $K$ in terms of their 
$P$ and $S$ components. Terms such as 
$\Pi_P K_{SP}$ and $A_{PP}K_{SP}$
vanish. By \eqref{APSmap} and \eqref{Kcond2}, the terms $K_{PS}A_{SP}$ 
and $K_{SP}A_{PS}$ are bounded on $L^2(\R^n)$, hence are absorbed into $R$. Using
\eqref{Kcond1}, the conditions \eqref{conj1}-\eqref{conj1'} thus reduce to the following conditions
\begin{align}
K_{PS}A_{SS}-A_{PP}K_{PS}+A_{PS}&= R\,,\label{conj2}\\
\rule{0pt}{16pt}
K_{SP}A_{PP}-A_{SS}K_{SP}+A_{SP}&= R\,,\label{conj2'}
\end{align}
meaning precisely that the left hand sides are bounded maps on $L^2(\R^n)$.
Condition \eqref{conj2'} is equivalent to \eqref{conj2} by taking adjoints, so we
will establish \eqref{conj2}.

Consider first the case of smooth $\mu$ and $\l$.
Since $K_{PS}$ is of order $-1$, and $A_{PP}$ and $A_{SS}$ of order 2, the
equality \eqref{conj2} reduces to saying that the top (i.e.~first) order terms in the 
symbol expansion for the left hand side cancel out.
By \eqref{princsymbA}, this implies that
\begin{equation}\label{KPSsymb}
\sigma_{-1}(K_{PS})=\frac{-\sigma_1(A_{PS})}{\bigl(\mu(x)+\l(x)\bigr)|\xi|^{2}}\,.
\end{equation}

We will also reduce \eqref{conj2} to a principal symbol calculation for the
case of $\mu$ and $\l$ of regularity $C^{1,1}$. Care needs to be taken,
however, since by \eqref{APSform} the symbol of $A_{PS}$, hence $K_{PS}$, is only 
Lipschitz in $x$, whereas $A_{PP}$ and $A_{SS}$ are of second order. 
That we can reduce matters to a symbol calculation
depends critically on the fact that $K_{PS}$
has a particular structure, which in turn depends on the 
special form of $A_{PS}$.

We consider the operator $A_{PS}$ more closely. Since $\div\,\Pi_S=0$, by \eqref{Aform}
we can write
\begin{align*}
(A_{PS}w)_i&=\sum_{j,k=1}^n\Delta^{-1}\partial_i\partial_j\partial_k
\bigl(\,\mu\,\bigl(\partial_k(\Pi_Sw)_j+\partial_j(\Pi_Sw)_k\bigr)\,\bigr)\\
&=2\sum_{j,k=1}^n\Delta^{-1}\partial_i\partial_j
\bigl(\,(\partial_k\mu)\,\partial_j(\Pi_Sw)_k\bigr)\\
&=2\,\sum_{k=1}^n\partial_i
\bigl(\,(\partial_k\mu)\,(\Pi_Sw)_k\bigr)
-2\sum_{j,k=1}^n\Delta^{-1}\partial_i\partial_j
\bigl(\,(\partial_j\partial_k\mu)\,(\Pi_Sw)_k\bigr)
\end{align*}
In going from the first to the second line we used that $\sum_k\partial_k(\Pi_Sw)_k=0$, 
together with symmetry of the expression in $j,k$. The second term in the third row
is a bounded operator on $L^2$; hence we may write
$$
A_{PS}w=2\,\nabla\bigl(\,(\grad \mu)\cdot\Pi_Sw\bigr)+R\,.
$$
Motivated by \eqref{KPSsymb}, we now define
$$
K_{PS}w=2\,\Phi(D)\Delta^{-1}\nabla\bigl(\,(\mu+\l)^{-1}(\nabla\mu)\cdot\Pi_Sw\bigr)\,,
$$
and set $K_{SP}=-K_{PS}^*$ so that \eqref{Kcond1} holds.
Here, $\Phi(\xi)$ is a smooth cutoff to $|\xi|\ge M$, 
$$
\Phi(\xi)=\begin{cases}1\,,&\quad |\xi| \ge 2M\,,\\ 0\,,&\quad |\xi|\le M\,,\end{cases}
$$
with $M$ depending on the
$C^{1,1}$ norms of $\mu$ and $\l$. Condition 
\eqref{Kcond2} holds since $(\l+\mu)^{-1}\nabla \mu$ is a multiplier on $H^s$ for
$-1\le s\le 1$. Furthermore, by taking $M$ large, we will have
\begin{equation}\label{Kcond3}
\|Kw\|_{H^s}\le \tfrac 12  \|w\|_{H^s}\,,\qquad -1\le s\le 1\,.
\end{equation}
Finally, we note that since $K_{PS}w$ is the gradient of a function, it follows that
$$
K_{PS}=\Pi_PK_{PS}=K_{PS}\Pi_S\,.
$$

To evaluate the composition $A_{PP}K_{PS}$, we observe that we may write
\begin{equation}\label{APP}
A_{PP}=(2\mu+\l)\,\Delta\Pi_P+\sum_{k=1}^nR_k\,\partial_k\,,
\end{equation}
where 
$$
R_k\,:\,H^s(\R^n)\rightarrow H^s(\R^n)\,,\qquad -1\le s\le 1\,.
$$
This follows from the fact that
the commutator $[\Pi_P\partial_k,\mu]$ is bounded on $H^s(\R^n)$ over this range, which is a
consequence of the Calder\'on commutator theorem; see Calder\'on (1965).

Since $\partial_k\circ K_{PS}$ is bounded on $L^2$,
it follows that
\begin{align*}
A_{PP}K_{PS}w&=2\,(2\mu+\l)\Phi(D)\nabla\bigl(\,(\mu+\l)^{-1}(\nabla\mu)\cdot\Pi_Sw\bigr)+Rw\\
&\rule{0pt}{14pt}
=2\,\Phi(D)\nabla\bigl(\,(2\mu+\l)\,(\mu+\l)^{-1}(\nabla\mu)\cdot\Pi_Sw\bigr)+Rw\,.
\end{align*}

To evaluate the composition $K_{PS}A_{SS}$, we first observe that we can write
$$
A_{SS}w=\Pi_S\Delta\bigl(\mu\,\Pi_S w\bigr)+A_1w\,,
$$
where $A_1\,:\,L^2\rightarrow H^{-1}$.
To see this, we note that $\div\,\Pi_S\equiv 0$, so the $\lambda$ terms in $A$ do
not contribute to $A_{SS}$. For the same reason we can write
$$
\Pi_S\sum_k\partial_k\bigl(\,\mu\,\partial_i(\Pi_Sw)_k\bigr)=
\Pi_S\sum_k(\partial_k\mu)\,\partial_i(\Pi_Sw)_k
$$
which maps $L^2\rightarrow H^{-1}$. Finally, we write
$$
\Pi_S\sum_k\partial_k\bigl(\,\mu\,\partial_k(\Pi_Sw)_i\bigr)=
\Pi_S\Delta(\,\mu\,\Pi_Sw)_i-
\Pi_S\sum_k(\partial_k\mu)\,\partial_k(\Pi_Sw)_i\,,
$$
and observe that the last term maps $L^2\rightarrow H^{-1}$.

By \eqref{Kcond2}, the composition $K_{PS}A_1$ is bounded on $L^2$, hence we
are reduced to showing that
\begin{equation}\label{KPSASScomp}
\Phi(D)\Delta^{-1}\nabla\bigl(\,(\mu+\l)^{-1}(\nabla\mu)\cdot(\Delta-\nabla\,\div)(\mu\,\Pi_Sw)\bigr)\\
=
\Phi(D)\nabla\bigl(\,(\mu+\l)^{-1}(\nabla\mu)\cdot(\mu\,\Pi_Sw)\bigr)+Rw\,,
\end{equation}
where we write $\Pi_S\Delta=\Delta-\nabla\,\div$, and understand that $\Delta$ acts 
on each component of a vector, i.e.~as a scalar operator.

The relation \eqref{KPSASScomp} is derived as a consequence of the identities
\begin{align}\label{ident1}
\Delta(V\!\cdot W)&= V\!\cdot\Delta W+2\,\div\,(W\!\cdot\nabla V)-W\!\cdot(\Delta V)
\\
\rule{0pt}{16pt}
\label{ident2}
V\!\cdot \nabla\,\div\,W&=\div\,(V\div\,W)-\div\,(W\div\, V)+W\!\cdot\nabla\,\div\,V
\end{align}
where in the term $\div\,(W\!\cdot\nabla V)$ the $\div$ pairs with $\nabla$.

Setting $V=(\mu+\l)^{-1}\nabla\mu$, and $W=\mu\,\Pi_Sw$, then \eqref{KPSASScomp} is equivalent
to
$$
\Phi(D)\Delta^{-1}\nabla\Bigl(\,(V\!\cdot(\Delta-\nabla\,\div)W-\Delta(V\!\cdot W)\,\Bigr)=Rw\,.
$$
By \eqref{ident1}-\eqref{ident2} we can write the left hand side as
$$
\Phi(D)\Delta^{-1}\nabla\Bigl(\,W\!\cdot(\Delta-\nabla\,\div)V-\div\,(V\div\,W)+\div(W\div V\,)
-2\,\div(W\!\cdot\nabla V)\,\Bigr)\,.
$$
The last two terms,
$$
\Phi(D)\Delta^{-1}\nabla\,\div\bigl(\,W\div\,V-2\,W\!\cdot\nabla V\bigr)
$$
are $L^2$ bounded operators applied to $W$, hence of the form $Rw$, since $V$ is Lipschitz.
Since $\div\,W=(\nabla\mu)\cdot\Pi_Sw$, the same holds for $\Phi(D)\Delta^{-1}\nabla\,\div(V\div\,W)$.
The remaining term is where we use the special structure, since it involves two derivatives of $V$. 
The key observation is that
$$
(\Delta-\grad\,\div)V=
(\Delta-\grad\,\div)\bigl(\,(\mu+\l)^{-1}\nabla\mu\,\bigr)\in L^\infty(\R^n)\,,
$$
since $\Delta\nabla\mu=\nabla\div\nabla\mu$, and hence the expression in fact involves only second
derivatives of $\mu$ and $\l$.
Consequently,
$$
w\rightarrow \Phi(D)\Delta^{-1}\nabla\Bigl(\,W\!\cdot(\Delta-\nabla\,\div)V\,\Bigr)\,:\,
L^2(\R^n)\rightarrow H^1(\R^n)\,.
$$

\newsection{The Parametrix and Decoupling of Modes}
We will realize the parametrix construction for \eqref{eq:CIVP} using
parametrices for the purely polarized equations, following the procedure introduced
in Smith (1998).
\begin{theorem}\label{parametrix}
There exist strongly continuous one parameter families of operators
\begin{equation*}
\begin{split}
\bc_P(t)&:H^s(\R^n)\rightarrow H^s(\R^n)\,,\qquad \qquad
\bs_P(t):H^{s-1}(\R^n)\rightarrow H^s(\R^n)\,,\\
\rule{0pt}{16pt}
\partial_t\bc_P(t)&:H^s(\R^n)\rightarrow H^{s-1}(\R^n)\,,\qquad \,
\partial_t\bs_P(t):H^{s-1}(\R^n)\rightarrow H^{s-1}(\R^n)\,,
\end{split}
\end{equation*}
for $0\le s \le 2$ and $t\in [-t_0,t_0]$, such that the solution to the Cauchy problem
$$
\bigl(\partial_t^2-A_{PP}\bigr)v_P=G_P\,,\qquad v_P|_{t=0}=f_P\,,
\quad \partial_t v_P|_{t=0}=g_P\,,
$$
where $\Pi_S f_P=\Pi_S g_P=\Pi_S G_P=0$, and
$$
f_P\in H^s(\R^n)\,,\qquad g_P\in H^{s-1}(\R^n)\,,\qquad 
G_P\in L^1\bigl([-t_0,t_0],H^{s-1}(\R^n)\bigr)\,,
$$
is given by
$$
v_P(t,\cdot\,)=\bc_P(t)f_P+\bs_P(t)g_P+\int_0^t s_P(t-s)\,G_P(s,\cdot\,)\,ds\,.
$$
The same result holds if $P$ is exchanged with $S$, and both results hold for all $t_0>0$.
\end{theorem}
\noindent{\it Proof.} Following pages 816--818 of Smith (1998)
we produce families of operators
$\rc(t)$ and $\rs(t)$, respectively bounded and regularizing of degree 1 on $H^s(\R^n)$ for all
$s$, with $\rc(0)=I$, $\partial_t\rc(0)=0$, $\rs(0)=0$, $\partial_t \rs(0)=I$,
such that for $0\le s\le 3$,
\begin{equation}\label{cp1}
\begin{split}
\bigl(\partial_t^2-(2\mu+\l)\Delta\bigr)\rc(t)=&\,T_0(t)\;:\,H^{s}(\R^n)\rightarrow H^{s-1}(\R^n)\,,\\
\rule{0pt}{16pt}
\bigl(\partial_t^2-(2\mu+\l)\Delta\bigr)\rs(t)=&\;T_1(t)\,:\,H^{s-1}(\R^n)\rightarrow H^{s-1}(\R^n)\,.
\end{split}
\end{equation}
We set $\rc_P(t)=\Pi_P\rc(t)\Pi_P$, $\rs_P(t)=\Pi_P\rs(t)\Pi_P$, and observe that, for
$0\le s\le 2$,
\begin{equation}\label{cp2}
\begin{split}
\bigl(\partial_t^2-A_{PP}\bigr)\rc_P(t)=&\,T_{P,0}(t)\;:\,H^{s}(\R^n)\rightarrow H^{s-1}(\R^n)\,,\\
\rule{0pt}{16pt}
\bigl(\partial_t^2-A_{PP}\bigr)\rs_P(t)=&\;T_{P,1}(t)\,:\,H^{s-1}(\R^n)\rightarrow H^{s-1}(\R^n)\,.
\end{split}
\end{equation}
To see this, we use the decomposition \eqref{APP} to see that
$$
A_{PP}\rc_P(t)-(2\mu+\l)\Delta\rc_P(t)\,:\,H^{s}(\R^n)\rightarrow H^{s-1}(\R^n)\,,\quad
0\le s\le 2\,.
$$
As in the discussion following \eqref{APP}, the commutator 
$$
[\Pi_P,(2\mu+\l)\Delta]\,: H^{s}(\R^n)\rightarrow H^{s-1}(\R^n)\,,\quad 0\le s\le 2\,,
$$
and \eqref{cp2} follows from \eqref{cp1}. The exact evolution operators $\bc_P$ and $\bs_P$
are obtained from $\rc_P$ and $\rs_P$ following pages 819--820 of Smith (1998),
using a convergent iteration. For example,
$$
\bc_P(t)=\rc_P(t)+\int_0^t\rs_P(t-s)\,T(s)\,ds\,,
$$
where $T(t)$ is the solution to
$$
T(t)+\int_0^tT_{P,1}(t-s)\,T(s)\,ds=-T_{P,0}(t)\,.
$$
We observe here that, since $\rc_P(t)$, $\rs_P(t)$ and $T_{P,0}(t)$ are orthogonal to $\Pi_S$ on
both the right and left, it follows that so is $\bc_P(t)$.

The bounds on $\partial_t\bc_P(t)$, etc., follow similarly, using the fact that
$\partial_t\rc(t)$ and $\partial_t\rs(t)$ respectively map $H^{s}(\R^n)\rightarrow H^{s-1}(\R^n)$ and
$H^{s-1}(\R^n)\rightarrow H^{s-1}(\R^n)$ for all $s$. $\qquad\square$

\bigskip

Suppose that $w$ is the solution to the linear elastic system \eqref{eq:CIVP}.
We pose
\begin{equation*}
v=(I+K)w\,,\qquad v_P=\Pi_Pv\,,\quad v_S=\Pi_Sv\,.
\end{equation*}
By \eqref{conj1} we can then write
\begin{align*}
\partial_t^2 v_P&=\Pi_P(I+K)A(x,D)w\\
{}&=A_{PP}(I+K)w+Rw\\
{}&=A_{PP}v_P+R(I+K)^{-1}v\,,
\end{align*}
where at the last stage we use \eqref{Kcond3} to invert $(I+K)$ on the range of Sobolev spaces
we will consider for $w$.
Similar considerations apply to $v_S$, hence we have the system
\begin{equation}\label{eq:decoupled}
\begin{cases}
\;\bigl(\partial_t^2 -A_{PP}\bigr)v_P=R_P v\,,\quad 
&v_P|_{t=0}=\tilde f_P\,,\quad\partial_tv_P|_{t=0}=\tilde g_P\,,\\
\rule{0pt}{18pt}
\;\bigl(\partial_t^2 -A_{SS}\bigr)v_S=R_S v\,,\quad 
&v_S|_{t=0}=\tilde f_S\,,\quad\partial_tv_S|_{t=0}=\tilde g_S\,,
\end{cases}
\end{equation}
where $R_P$ and $R_S$ 
are bounded on $L^2(\R^n)$. 
Here, we have set $\tilde f_P=\Pi_P(I+K)f$ and $\tilde g_P=\Pi_P(I+K)g$, 
and similarly for $\tilde f_S$, $\tilde g_S$.

The system \eqref{eq:decoupled} is equivalent to the integral system,
\begin{equation}\label{eq:integral}
\begin{cases}
\;\;v_P(t)=\bc_P(t)\tilde f_P+\bs_P(t)\tilde g_P+\displaystyle\int_0^t \bs_P(t-s) R_P v(s,\cdot\,)\,ds
\,,\\\rule{0pt}{25pt}
\;\;v_S(t)=\bc_S(t)\tilde f_S+\bs_S(t)\tilde g_S+\displaystyle\int_0^t \bs_S(t-s) R_S v(s,\cdot\,)\,ds
\,.
\end{cases}
\end{equation}
We now restrict attention to data $f\in H^s(\R^n)$, $g\in H^{s-1}(\R^n)$, with $0\le s\le 1$.
Then
$$
\bc_P(t)\tilde f_P+\bs_P(t)\tilde g_P\in C\bigl([-t_0,t_0],H^s(\R^n)\bigr)\,,
$$
and the same holds with $P$ replaced by $S$. For $0\le s\le 1$, then $s_P(t-s)R_P$ is
a bounded map on $H^s(\R^n)$, and similarly with $P$ replaced by $S$, so the Volterra type
system \eqref{eq:integral} can be
solved by iteration to yield a solution 
$$
v=v_P+v_S\in C\bigl([-t_0,t_0],H^s(\R^n)\bigr)\,.
$$
The solution takes the form
$$
v=\bc_P(t)\tilde f_P+\bs_P(t)\tilde g_P+\bc_S(t)\tilde f_S+\bs_S(t)\tilde g_S+\br(t)(\tilde f,\tilde g)
\,,
$$
where $\br(t):L^2(\R^n)\oplus H^{-1}(\R^n)\rightarrow H^1(\R^n)$
and $\partial_t\br(t):L^2(\R^n)\oplus H^{-1}(\R^n)\rightarrow L^2(\R^n)$, with both $\br(t)$
and $\partial_t\br(t)$ strongly continuous in $t$.

We next write $\tilde f=(I+K)f$, $\tilde g=(I+K)g$, and $w=(I+K)^{-1}v$.
By \eqref{Kcond2},
$$
(I+K)^{-1}-I\,:\,H^s(\R^n)\rightarrow H^{s+1}(\R^n)\quad\text{for}\quad -1\le s\le 1\,.
$$
Consequently,
$$
(I+K)^{-1}\bc_P(t)(I+K)-\bc_P(t)\,:\,L^2(\R^n)\rightarrow H^1(\R^n)\,.
$$
By similar results for the other terms, we have the following corollary
(where $\br(t)$ is of a different form than above).
\begin{theorem}\label{wdecouple}
Let $w(t,x)$ solve the Cauchy problem \eqref{eq:CIVP}, with data $f\in L^2(\R^n)$, 
$g\in H^{-1}(\R^n)$. Let $(f_P,g_P)=\Pi_P(f,g)$, and $(f_S,g_S)=\Pi_S(f,g)$. 
Then $w$ may be written as
$$
w(t,\cdot\,)=\bc_P(t)f_P+\bs_P(t)g_P+\bc_S(t)f_S+\bs_S(t)g_S+\br(t)(f,g)\,,
$$
where $\rule{0pt}{12pt}\mathbf{r}(t)$ is bounded
from $L^2(\R^n)\oplus H^{-1}(\R^n)$ to $H^1(\R^n)$, and
$\partial_t\mathbf{r}(t)$ is bounded
from $L^2(\R^n)\oplus H^{-1}(\R^n)$ to $L^2(\R^n)$.
Both $\br(t)$ and $\partial_t\br(t)$ are strongly continuous in $t$.
\end{theorem}

In particular, if the initial data for $w$ is a pure $P$ mode, in that
$f_S=g_S=0$, then
$$
\Pi_Sw=\Pi_S\br(t)(f,g)\in C\bigl([-t_0,t_0],H^1(\R^n)\bigr)\,,
$$
provided that $f_P\in L^2(\R^n)$ and $g_P\in H^{-1}(\R^n)$. This can be viewed
as saying that the $P\rightarrow S$ coupling is regularizing of order 1, and
the same holds for pure $S$ modes.

This decoupling can also be expressed in terms of mapping properties of the components
for the solution operator to \eqref{eq:CIVP}. 
An immediate result of Theorem \ref{wdecouple} is the following.
\begin{corollary}\label{cor:decoupled}
Let the solution to \eqref{eq:CIVP}
be written in the form
$$
w(t,\cdot\,)=\bc(t)f+\bs(t)g\,,
$$
and decompose $\bc=\bc_{PP}+\bc_{PS}+\bc_{SP}+\bc_{SS}$, 
$\bs=\bs_{PP}+\bs_{PS}+\bs_{SP}+\bs_{SS}$, where $\bc_{PP}=\Pi_P\bc\,\Pi_P$, etc.
Then with $\bc_P\,,\bc_S\,,\bs_P\,,\bs_S$ as in Theorem \ref{parametrix}, we have

\begin{align*}
\bc_{PP}(t)-\bc_P(t)
&:L^2(\R^n)\rightarrow H^1(\R^n)\,,\rule{0pt}{12pt}\\
\bc_{PS}(t)
&:L^2(\R^n)\rightarrow H^1(\R^n)\,,\rule{0pt}{12pt}\\
\rule{0pt}{23pt}
\bs_{PP}(t)-\bs_P(t)
&:H^{-1}(\R^n)\rightarrow H^1(\R^n)\,,\rule{0pt}{12pt}\\
\bs_{PS}(t)
&:H^{-1}(\R^n)\rightarrow H^1(\R^n)\,.\rule{0pt}{12pt}
\end{align*}
All results hold with $P$ and $S$
interchanged, and are strongly continuous in $t$.

Furthermore, in each case 
the time derivative of the operator is a bounded map from
the indicated domain into $L^2(\R^n)$. For example,
$$
\partial_t\bc_{PP}(t)-\partial_t\bc_P(t)
:L^2(\R^n)\rightarrow L^2(\R^n)\,,
$$
with strongly continuous dependence on $t$.
\end{corollary}
The corollary follows from Theorem \ref{wdecouple}
since all operators listed arise as compositions of $\mathbf{r}(t)$ with
the projections $\Pi_P$ and $\Pi_S$.

We remark that if the stronger regularity assumption is made, $\mu\,,\l\in C^{k,1}$ with
$k\ge 2$, then the range of Sobolev indexes $s$ on which $\br(t)$ is regularizing of order 1
increases, and Corollary \ref{cor:decoupled} thus applies to a range of $H^s(\R^n)$.
Alternatively, by adding successive terms to $K$, the system \eqref{eq:decoupled}
for $v=(I+K)w$
can be decoupled to higher order, that is, with $R_P$ and $R_S$ of negative order; 
see Taylor (1975)
for the case of smooth $\mu$, $\l$ where $R_P$ and $R_S$ are smoothing operators.
This requires adding lower order terms to $A_{PP}$ and $A_{SS}$.
The order of gain in
regularity in Corollary \ref{cor:decoupled} cannot be increased, however, since it
is limited by the term of order $-1$ in $K$.

\bigskip
\bigskip

\noindent{\Large\bf Acknowledgments}

\medskip

The research of Brytik, de Hoop, and Uhlmann was supported in part
under NSF CMG grant EAR-0724644; Uhlmann was also partly supported
by a Chancellor Professorship at UC Berkeley and a Senior Clay Award.
The research of Smith was supported under NSF grant DMS-0654415.

\bigskip
\bigskip

\noindent{\Large\bf References}

\medskip


\noindent
Calder\'on, A.~(1965).
{Commutators of singular integral operators}. 
{\it Proc.~Nat.~Acad.~Sci.~U.S.A.} 53:1092--1099.

\smallskip\noindent
Dencker, N.~(1982).
{On the propagation of polarization sets for systems of real principal type}.
{\it J.~Functional Analysis} 46:351--372.

\smallskip\noindent
Langston, C.A.~(1979).
{Structure under {M}ount {R}ainier, {W}ashington, inferred from teleseismic body waves}.
{\it J.~Geophys.~Res.} 84:4749--4762.

\smallskip\noindent
Lin, C.L., Nakamura, G., Uhlmann, G., Wang, J.N.~(2010).
{Quantitative unique continuation for the Lam\'e system with less regular coefficients}.
Preprint, arXiv:1005.3382v1.

\smallskip\noindent
Smith, H.F.~(1998).
{A parametrix construction for wave equations with {$C^{1,1}$} coefficients}.
{\it Ann.~Inst.~Fourier, Grenoble} 4:797--835.

\smallskip\noindent
Stolk, C.~(2000).
{Microlocal analysis of a seismic linearized inverse problem}.
{\it Wave Motion} 32:267--290.

\smallskip\noindent
Stolk, C., de Hoop, M.V.~(2002).
{Microlocal analysis of seismic inverse scattering in anisotropic elastic media}.
{\it Comm.~Pure Appl.~Math.} 55:261--301.

\smallskip\noindent
Stolk, C., de Hoop, M.V.~(2006).
{Seismic inverse scattering in the downward continuation approach}.
{\it Wave Motion} 43:579--598.

\smallskip\noindent
Taylor, M.~(1975).
{Reflection of singularities of solutions to systems of differential equations}.  
{\it Comm.~Pure Appl.~Math.} 28:457--478.

\smallskip\noindent
Wang, J.N.~(1998).
{On the relation between singularities of coefficients and singularities of 
reflected waves in the Lam\'e system}. 
{\it Inverse Problems} 14:733--743.

\smallskip\noindent
Vinnik, L.~(1977).
{Detection of waves converted from {P} to {SV} in the mantle}.
{\it Phys. Earth Planet. Inter.} 15:39--45.


\end{document}